\documentclass[english]{svjour3}
\usepackage[latin9]{inputenc}
\usepackage{geometry}
\usepackage{babel}
\usepackage{url}
\usepackage{multirow}
\usepackage{amsmath}
\usepackage{graphicx}
\usepackage{esint}
\usepackage[numbers]{natbib}
\usepackage[unicode=true,pdfusetitle,
 bookmarks=true,bookmarksnumbered=true,bookmarksopen=true,bookmarksopenlevel=2,
 breaklinks=false,pdfborder={0 0 1},backref=false,colorlinks=false]
 {hyperref}
\usepackage{breakurl}

\makeatletter

\providecommand{\tabularnewline}{\\}

\numberwithin{equation}{section}
\numberwithin{figure}{section}
  \newenvironment{svmultproof}{\begin{proof}}{\qed\end{proof}}

\usepackage{babel}

\usepackage{babel}

\usepackage{amsfonts}
\usepackage{graphicx}
 \usepackage{layout}
\usepackage{caption}

\makeatother

\begin{document}

\title{Banded operational matrices for Bernstein polynomials and application
to the fractional advection-dispersion equation}

\titlerunning{Banded operational matrices for Bernstein polynomials and application
...}

\author{M. Jani \and E. Babolian \and S. Javadi \and D. Bhatta}

\institute{M. Jani \and E. Babolian \and S. Javadi \at Department of Mathematics,
Faculty of Mathematical Sciences and Computer, Kharazmi University,
Tehran, Iran\\
\email{mostafa.jani@gmail.com babolian@khu.ac.ir javadi@khu.ac.ir}\\
\and D. Bhatta \at School of Mathematical \& Statistical Sciences,
The University of Texas Rio Grande Valley, 1201 West University Drive,
Edinburg, TX, USA\\
\email{dambaru.bhatta@utrgv.edu} \\
}
\maketitle
\begin{abstract}
In the papers dealing with derivation and applications of operational
matrices of Bernstein polynomials, a basis transformation, commonly
a transformation to power basis, is used\@. The main disadvantage
of this method is that the transformation may be ill-conditioned\@.
Moreover, when applied to the numerical simulation of a functional
differential equation, it leads to dense operational matrices and
so a dense coefficient matrix is obtained\@. In this paper, we present
a new property for Bernstein polynomials\@. Using this property,
we build exact banded operational matrices for derivatives of Bernstein
polynomials\@. Next, as an application, we propose a new numerical
method based on a Petrov-Galerkin variational formulation and the
new operational matrices utilizing the dual Bernstein basis for the
time-fractional advection-dispersion equation\@. We show that the
proposed method leads to a narrow-banded linear system and so less
computational effort is required to obtain the desired accuracy for
the approximate solution\@. We also obtain the error estimation for
the method\@. Some numerical examples are provided to demonstrate
the efficiency of the method and to support the theoretical claims\@.

\end{abstract}

\section{\label{sec:Intro}Introduction}

Bernstein polynomials play an important role in computer aided geometric
design \citep{farin2002}\@. Moreover, these polynomials provide
a tool in the numerical simulation of differential, integral and integro-differential
equations; see e.g., \citep{malek2012,maleknejad2011,Youzbasi} and
the references therein\@. In recent decades, many authors discovered
some new analytic properties and also some applications for these
polynomials\@. For example, Cheng \citep{cheng1983} derived the
rate of convergence of these polynomials for a certain class of functions\@.
Farouki \citep{farouki2-1996} showed that among nonnegative bases
on a given interval, the Bernstein polynomials basis is an optimal
stable basis\@. Delgado et al. \citep{delgado2009} proved that the
collocation matrices of the Bernstein basis are the best conditioned
among all the collocation matrices of nonnegative totally positive
bases on a given interval\@.

On the other hand, a basis transformation, commonly to the power basis
$\{1,x,\dots,x^{N}\}$, is used in order to derive the operational
matrices for the derivatives and integrals of Bernstein polynomials\@.
The developed numerical method for differential and integral equations
based on those matrices leads to a linear system whose coefficient
matrix is neither banded nor sparse; see, for instance, \citep{malek2012,pandey,parand,saadatmandi2014,Yousefi2010,Yousefi2011}\@.
It is also worth to note that the explicit conversion between the
Bernstein basis and the power basis is exponentially ill-conditioned
\citep{Farouki1999}\@.

In this paper, we directly derive narrow-banded operational matrices
for the derivatives of Bernstein polynomials without using any basis
transformation and we show that it leads to less computational cost
and less round-off errors\@. Then, as an application, we propose
a numerical scheme for the time fractional advection-dispersion equation
(FAD)\@. In the proposed method, we use the Bernstein polynomials
as the trial functions and the dual Bernstein polynomials as the test
functions for the Petrov-Galerkin variational formulation\@. We then
show that under some reasonable assumptions, the derived linear system
has a unique solution\@.

Due to various applications of fractional advection-dispersion equation
in the mathematical modeling of some important physical problems,
many classical numerical methods for the partial differential equations
have been developed to handle this problem\@. For example, Gao et
al. \citep{Gao} proposed a numerical scheme based on a compact finite
difference method\@. Shirzadi et al. \citep{Shirzadi} implemented
the meshless methods\@. Jiang et al. \citep{Jiang} developed a finite
element method\@. Stokes et al. \citep{Stok} used a mapping from
the associated Brownian counterpart for the numerical solution of
the time-fractional advection-dispersion equation\@.

This paper is concerned with providing a new numerical method for
the time-fractional advection-dispersion equation given by
\begin{eqnarray}
 &  & {D_{t}^{\alpha}}u\left(x,t\right)=\kappa_{1}\frac{\partial^{2}}{\partial x^{2}}u\left(x,t\right)-\kappa_{2}\frac{\partial}{\partial x}u\left(x,t\right)+S(x,t),\quad\left(x,t\right)\in\Omega\times\left(0,\infty\right),\label{eq:main}
\end{eqnarray}
where $u$ is the unknown concentration, $\kappa_{1}$ and $\kappa_{2}$
are the dispersion and advection coefficients, $\alpha\in\left(0,1\right)$
is the temporal fractional order, $S$ is the known source term, the
operator $D_{t}^{\alpha}$ represents the Caputo fractional derivative
and $\Omega=\left(0,L\right).$
\begin{definition}
\label{CaputoDef}\citep{Diethelm} The Caputo fractional derivative
of order $\alpha>0$ is defined by 
\begin{eqnarray*}
 &  & {D^{\alpha}}f(x)=\frac{1}{\Gamma(m-\alpha)}\int_{0}^{x}{\frac{1}{(x-s)^{\alpha-m+1}}\frac{d^{m}f(s)}{ds^{m}}ds},\quad x>0,
\end{eqnarray*}
where $m:=\left\lceil \alpha\right\rceil $\@.
\end{definition}
So the Caputo temporal fractional derivative of order $\alpha$ is
written as
\begin{eqnarray}
 &  & {D_{t}^{\alpha}}u(x,t)=\frac{1}{\Gamma(m-\alpha)}\int_{0}^{t}{\frac{1}{(t-s)^{\alpha-m+1}}\frac{\partial^{m}u(x,s)}{\partial s^{m}}ds}.\label{eq:4 CaputoFracDef}
\end{eqnarray}
We consider the equation (\ref{eq:main}) with the following initial
and boundary conditions 
\begin{eqnarray}
 &  & u\left(x,0\right)=g\left(x\right),\quad x\in\left[0,L\right],\label{IV}\\
 &  & u\left(0,t\right)=0,\quad u\left(L,t\right)=0,\quad t>0.\label{BVs}
\end{eqnarray}
This paper is organized as follows\@. In Section \ref{sec:PreliFracBern},
we introduce some background material of Bernstein polynomials\@.
Section \ref{sec:BernDerivatives=000026Matrices} is devoted to some
new relations for Bernstein  basis and the derivation of new banded
operational matrices\@. In Section \ref{sec:MainFormulationAndApplication},
we use the new operational matrices and the Petrov-Galerkin variational
formulation to develop a numerical simulation for the time-fractional
advection-dispersion equation\@. The error analysis for the method
is discussed in Section \ref{sec:Error}\@. Some numerical examples
are carried out in Section \ref{sec:Error}\@. The paper ends with
some concluding remarks in Section \ref{sec:Con}\@.

\section{\label{sec:PreliFracBern}Preliminaries}

We first provide the definition and some preliminary results for Bernstein
polynomials\@.

The Bernstein polynomials of degree $N\geq0$ on $[a,b]$ are defined
by
\begin{eqnarray}
 &  & B_{i,N}(x)=\frac{\dbinom{N}{i}(x-a)^{i}(b-x)^{N-i}}{(b-a)^{N}},\qquad0\leq i\leq N.\label{eq:1 BernDef}
\end{eqnarray}
We adopt the convention $B_{i,N}(x)\equiv0$ for $i<0$ and $i>N$\@.
The set $B_{N}=\left\{ B_{i,N}(x):\,i=0,\dots,N\right\} $ is a basis
for $\mathbb{P}_{N}$, the space of polynomials of degree not exceeding
$N$\@. 

The base functions satisfy the endpoint interpolation property, i.e.
\begin{eqnarray}
 &  & B_{i,N}(a)=\delta_{i,0},\quad B_{i,N}(b)=\delta_{i,N},\quad0\leq i\leq N,\,N\geq0.\label{eq:BernBound}
\end{eqnarray}
It is useful in the numerical formulation of the problem when imposing
boundary conditions\@. Also, for $N\geq1,$ the basis has the following
degree elevation property:
\begin{eqnarray*}
 &  & B_{i,N-1}(x)=\frac{1}{N}\left[\left(N-i\right)B_{i,N}(x)+\left(i+1\right)B_{i+1,N}(x)\right],\quad0\leq i\leq N.
\end{eqnarray*}
This relation may be recursively used to get the following general
formula\@.
\begin{lemma}
{[}\citealp{farouki1988}, Relation (26){]} Let $i,\,j$ and $N$
be nonnegative integers, $j\leq N$ and $i\leq N-j$\@. Then, 
\begin{eqnarray}
 &  & B_{i,N-j}(x)=\binom{N-j}{i}\sum_{r=i}^{j+i}\frac{\binom{j}{r-i}}{\binom{N}{r}}B_{r,N}(x).\label{eq:DegreeElevation}
\end{eqnarray}
\end{lemma}
The derivatives of the Bernstein polynomials satisfy the following
recurrence relation {[}\citealp{farouki1988}, Relation (54){]} 
\begin{eqnarray*}
 &  & B_{i,N}^{\prime}(x)=\frac{N}{b-a}\left(B_{i-1,N-1}(x)-B_{i,N-1}(x)\right),\quad0\leq i\leq N.
\end{eqnarray*}
Using the Leibniz's rule, we derive the following result (see {[}\citealp{doha2-2011},
Theorem 3.1.{]} for a similar result on the unit interval)\@.
\begin{lemma}
\label{Lem:HighorderDerChangeBasis}Let $N$ and $p$ be nonnegative
integers and $p\leq N$\@. Then, for $0\leq i\leq N,$ 
\begin{eqnarray}
 &  & B_{i,N}^{(p)}(x)=c_{p,N}\sum_{k=\max{(0,i+p-N)}}^{\min{(i,p)}}(-1)^{k}\dbinom{p}{k}B_{i-k,N-p}(x),\label{eq:HighDer.DegChange}
\end{eqnarray}
where
\begin{eqnarray}
 &  & c_{p,N}=\frac{\left(-1\right)^{p}N!}{(b-a)^{p}(N-p)!}.\label{Symbol_c_=00005Bp,N=00005D}
\end{eqnarray}
\end{lemma}
\begin{theorem}
{[}\citealp{juttler}, Theorem 3.{]} The elements of the dual basis
$B_{N}^{\star}=\left\{ B_{i,N}^{\star}(x):\,i=0,\dots,N\right\} $
associated with the Bernstein basis $B_{N}$ on $[a,b]$ are given
by 
\begin{eqnarray}
 &  & B_{i,N}^{\star}(x)=\sum_{j=0}^{N}d_{i,j}B_{j,N}(x),\quad0\leq i\leq N,\label{eq:Dual}
\end{eqnarray}
with the coefficients 
\begin{eqnarray*}
 &  & d_{i,j}=\frac{(-1)^{i+j}}{(b-a)\binom{N}{i}\binom{N}{j}}\sum_{r=0}^{min(i,j)}(2r+1)\binom{N+r+1}{N-i}\binom{N-r}{N-i}\binom{N+r+1}{N-j}\binom{N-r}{N-j},
\end{eqnarray*}
for $i,j=0,1,\dotsc,N$\@. Two sets $B_{N}$ and $B_{N}^{\star}$
form a biorthogonal system, i.e.
\begin{eqnarray}
 &  & \int_{a}^{b}B_{i,N}(x)B_{j,N}^{\star}(x)dx=\delta_{ij},\label{eq2.11}
\end{eqnarray}
for $i,j=0,1,\dots,N$\@.
\end{theorem}

\section{\label{sec:BernDerivatives=000026Matrices}Banded operational matrices
for derivatives of Bernstein polynomials }

The existing operational matrices for Bernstein basis and the applications
are based on a basis transformation, commonly from Bernstein to power
basis (see, for instance,\citep{malek2012,pandey,parand,saadatmandi2014,Yousefi2010,Yousefi2011}).
The transformation may be ill-conditioned, and also it does not lead
to banded linear systems \citep{Farouki1999}\@. In this section,
we first obtain a new representation for derivatives of Bernstein
polynomials and then present banded operational matrices\@.

The representation (\ref{eq:HighDer.DegChange}) for $p$th order
derivative of Bernstein polynomials changes the basis order from $N$
to $N-p$\@. For the Galerkin and Petrov-Galerkin formulations of
the equation (\ref{eq:main}), we need to express the derivatives
in terms of the same basis functions\@. So we state the following
theorem\@.
\begin{theorem}
\label{Theorem.Main}Let $N$ be any nonnegative integer and $0\leq i\leq N$\@.
Then, for any nonnegative integer $p\leq N,$ we have the following
``at most $(2p+1)$ term'' relation
\begin{eqnarray}
 &  & B_{i,N}^{(p)}(x)=c_{p,N}\sum_{j=\max(0,i-p)}^{\min(N,i+p)}\omega_{i,j}B_{j,N}(x),\quad0\leq i\leq N,\label{eq3.4-1-1}
\end{eqnarray}
where $c_{p,N}$ is defined as (\ref{Symbol_c_=00005Bp,N=00005D})
and the coefficients are as follows:
\begin{eqnarray}
 &  & \omega_{i,j}=\sum_{k=\max(0,i-j)}^{\min(p,i-j+p)}{\mu_{i,k,j}},\nonumber \\
 &  & \mu_{i,k,j}=(-1)^{k}\frac{\binom{p}{k}\binom{N-p}{i-k}\binom{p}{j-i+k}}{\binom{N}{j}}.\label{eq:muFormula}
\end{eqnarray}
\end{theorem}
\begin{svmultproof}
By using (\ref{eq:DegreeElevation}) and (\ref{eq:HighDer.DegChange}),
we obtain the following relation
\begin{eqnarray*}
B_{i,N}^{(p)}(x) & = & c_{p,N}\sum_{k=\max{(0,i+p-N)}}^{\min{(i,p)}}{(-1)^{k}\binom{p}{k}\binom{N-p}{i-k}\sum_{j=i-k}^{i-k+p}{\frac{\binom{p}{j-i+k}}{\binom{N}{j}}B_{j,N}(x)}}\\
 & = & c_{p,N}\sum_{k=m}^{n}{\left(\sum_{j=i-k}^{i-k+p}{\mu_{i,k,j}B_{j,N}(x)}\right)},
\end{eqnarray*}
where $m=\max{(0,i+p-N)}$ and $n=\min{(i,p)}$\@. The summation
limits are $m\leq k\leq n$ and $i-k\leq j\leq i-k+p$, changing the
order of summation as $i-n\leq j\leq i-m+p$ and $i-j\leq k\leq i-j+p$
and noting that $i+p-m=\min{(N,i+p)}$ and $i-n=\max{(0,i-p)},$ we
get 
\begin{eqnarray*}
B_{i,N}^{(p)}(x) & = & c_{p,N}\sum_{j=\max{(0,i-p)}}^{\min{(N,i+p)}}{\left(\sum_{k=i-j}^{i-j+p}{\mu_{i,k,j}}\right)B_{j,N}(x)}.
\end{eqnarray*}
Due to the fact that $\mu_{i,k,j}=0$ for $k<0$ and $k>p$, by removing
the zero terms, we obtain:
\begin{eqnarray*}
B_{i,N}^{(p)}(x) & = & c_{p,N}\sum_{j=\max{(0,i-p)}}^{\min{(N,i+p)}}{\left(\sum_{k=\max(0,i-j)}^{\min(p,i-j+p)}{\mu_{i,k,j}}\right)B_{j,N}(x)}\\
 & = & c_{p,N}\sum_{j=\max(0,i-p)}^{\min(N,i+p)}\omega_{i,j}B_{j,N}(x).
\end{eqnarray*}
This ends the proof\@.
\end{svmultproof}

\begin{remark}
For $p=1$ and $p=2$, the relation (\ref{eq3.4-1-1}) is written
as
\begin{eqnarray}
 &  & B_{i,N}^{\prime}(x)=\frac{1}{b-a}((N-i+1)B_{i-1,N}(x)-(N-2i)B_{i,N}(x)-(i+1)B_{i+1,N}(x)),\label{eq:firstDerCoef}\\
 &  & B_{i,N}^{\prime\prime}(x)=\frac{1}{\left(b-a\right)^{2}}((N-i+2)(N-i+1)B_{i-2,N}(x)-2(N-i+1)(N-2i+1)B_{i-1,N}(x)\,\label{eq:SecDerCoeff}\\
 &  & \qquad+(N^{2}-6Ni+6i^{2}-N)B_{i,N}(x)+2(i+1)(N-2i-1)B_{i+1,N}(x)+(i+2)(i+1)B_{i+2,N}(x)),\nonumber 
\end{eqnarray}
for $i=0,1,\dots,N$\@. 
\end{remark}
Let $\phi=[B_{i,N}(x):\,i=0,1,...,N]^{T}$\@. By using (\ref{eq3.4-1-1}),
we have 
\begin{eqnarray*}
\frac{d^{p}}{dx^{p}}\phi & = & \mathbf{D}_{p}\phi,\quad p\geq1,
\end{eqnarray*}
where $\mathbf{D}_{p}$ is a $(2p+1)$-diagonal, or $(p,p)$-banded,
matrix expressed as 
\[
\mathbf{D}_{p}=c_{p,N}\left[\begin{array}{ccccccccccc}
\omega_{0,0} &  & \cdots &  & \omega_{0,p}\\
\omega_{1,0} &  &  & \cdots &  & \omega_{1,p+1}\\
\vdots\\
\omega_{p,0} &  &  &  &  &  & \omega_{p,2p}\\
 & \omega_{p+1,1} &  &  &  & \cdots &  & \omega_{p+1,2p+1}\\
 &  & \omega_{p+2,2} &  &  &  & \cdots &  & \omega_{p+2,2p+2}\\
 &  &  & \ddots &  &  &  &  & \ddots\\
\\
\\
\\
 &  &  &  &  &  &  &  & \omega_{N,N-p} & \cdots & \omega_{N,N}
\end{array}\right].
\]
The above matrix, written according to Theorem \ref{Theorem.Main},
avoids matrix multiplications $\mathbf{D}_{p}=\mathbf{D}_{1}^{p}$\@.
\begin{remark}
We especially note that $\mathbf{D}_{0}=\mathbf{I}$, the identity
matrix, $\mathbf{D}_{1}$ is the following tridiagonal matrix:
\[
\mathbf{D}_{1}=\frac{1}{b-a}\left[\begin{array}{cccccc}
-N & -1 & 0 &  & \cdots & 0\\
N & 2-N & -2 & 0 &  & \vdots\\
0 & N-1 & 4-N & -3 & 0\\
 & 0 & \ddots & \ddots & \ddots & 0\\
\vdots &  & 0 & 2 & N-2 & -N\\
0 & \cdots &  & 0 & 1 & N
\end{array}\right],
\]
and $\mathbf{D}_{2}$ is a pentadiagonal matrix with elements given
by
\begin{eqnarray*}
\left(\mathbf{D}_{2}\right)_{i,j} & = & \frac{1}{\left(b-a\right)^{2}}\times\begin{cases}
\begin{array}{cc}
(N-i+2)(N-i+1), & i-j=-2,\\
-2(N-i+1)(N-2i+1), & i-j=-1,\\
N^{2}-6Ni+6i^{2}-N, & i-j=0,\\
2(i+1)(N-2i-1), & i-j=1,\\
(i+2)(i+1), & i-j=2,\\
0, & \left|i-j\right|>2.
\end{array} & i,j=0,\ldots.,N.\end{cases}
\end{eqnarray*}
\end{remark}

\begin{remark}
\label{Remark}For $\mathbf{D}_{1}$, it is seen that sum of the elements
in each column is zero\@. This holds for $\mathbf{D}_{p}$, $p\geq1$.
To see this, let $0\leq j\leq N$ and $p>1$\@. Then,
\begin{eqnarray*}
\sum_{i=0}^{N}{(\mathbf{D}_{p})_{i,j}} & = & \sum_{i=0}^{N}{(\mathbf{D}_{1}\mathbf{D}_{p-1})_{i,j}}=\sum_{i=0}^{N}{\sum_{k=0}^{N}{(\mathbf{D}_{1})_{i,k}(\mathbf{D}_{p-1})_{k,j}}}\\
 & = & \sum_{k=0}^{N}{(\mathbf{D}_{p-1})_{k,j}\sum_{i=0}^{N}{(\mathbf{D}_{1})_{i,k}}}=0.
\end{eqnarray*}
\end{remark}
Here we present some other interesting features of $\mathbf{D}_{p}$.
\begin{lemma}
\label{LemmaNilpotent}Let $p$ be a positive integer\@. Then $\mathbf{D}_{p}$
is a nilpotent matrix, i.e., $\mathbf{D}_{p}^{k}=0$ for some positive
integer $k$\@. In fact the smallest such $k$ is $k=\left\lceil \frac{N+1}{p}\right\rceil $
where $\left\lceil \cdot\right\rceil $ represents the ceiling function,
i.e., the smallest following integer\@.
\end{lemma}
\begin{svmultproof}
Let $i\leq j$. Then,
\begin{eqnarray*}
(\mathbf{D}_{1}^{N})_{i,j} & = & (\mathbf{D}_{N})_{i,j}=c_{N,N}\sum_{k=\max(0,i-j)}^{\min(N,i-j+N)}{\mu_{i,k,j}}\\
 & = & c_{N,N}\sum_{k=0}^{N+i-j}{(-1)^{k}\frac{\binom{N}{k}\binom{0}{i-k}\binom{N}{j-i+k}}{\binom{N}{j}}}\\
 & = & c_{N,N}(-1)^{i}\binom{N}{i}.
\end{eqnarray*}
Similarly, for $i>j$ we obtain $(\mathbf{D}_{1}^{N})_{i,j}=c_{N,N}(-1)^{i}\binom{N}{i}$,
and so
\begin{eqnarray}
(\mathbf{D}_{1}^{N})_{i,j} & = & c_{N,N}(-1)^{i}\binom{N}{i},\quad0\leq i,j\leq N.\label{eq:D1^N}
\end{eqnarray}
This means that for each row of $\mathbf{D}_{1}^{N}$, all the entries
are the same, and we obtain 
\begin{eqnarray*}
(\mathbf{D}_{1}^{N+1})_{i,j} & = & \sum_{k=0}^{N}{(\mathbf{D}_{1}^{N})_{i,k}(\mathbf{D}_{1})_{k,j}}\\
 & = & c_{N,N}(-1)^{i}\binom{N}{i}\sum_{k=0}^{N}{(\mathbf{D}_{1})_{k,j}},
\end{eqnarray*}
for $0\leq i,j\leq N+1$\@. From Remark \ref{Remark}, the last summation
is zero and so $\mathbf{D}_{1}^{N+1}=0$\@. It means that $\mathbf{D}_{1}$
is nilpotent\@. Since $\mathbf{D}_{1}^{N}\neq0$ by (\ref{eq:D1^N}),
$\mathbf{D}_{1}$ is a nilpotent matrix with index $N+1$\@. Let
$k=\left\lceil \frac{N+1}{p}\right\rceil $\@. Using $kp\geq N+1$
and the fact that $\mathbf{D}_{p}=\mathbf{D}_{1}^{p}$, we obtain
\begin{eqnarray*}
\mathbf{D}_{p}^{k} & = & \mathbf{D}_{1}^{kp}=0.
\end{eqnarray*}
On the other hand, since $(k-1)p<N+1$ so $\mathbf{D}_{p}^{k-1}=\mathbf{D}_{1}^{(k-1)p}\neq0$
and so $\mathbf{D}_{p}$ is a nilpotent matrix with index $k$\@.
This completes the proof\@.
\end{svmultproof}

The following corollary is obtained from the fact that $\mathbf{D}_{p}$
is nilpotent (see \citep{Serre}, for more properties of nilpotent
matrices)\@.
\begin{corollary}
Let $p\geq1$ and $0\leq i,j\leq N$, then

(a) All the eigenvalues of $\mathbf{D}_{p}$ are zero, 

(b) $trace(\mathbf{D}_{p})=0,$ 

(c) $\mathbf{I}-c\mathbf{D}_{p}$ is nonsingular for any scalar $c\in\mathbb{C}$
and 
\begin{eqnarray*}
(\mathbf{I}-c\mathbf{D}_{p})^{-1} & = & \mathbf{I}-c\mathbf{D}_{1}^{p}+c^{2}\mathbf{D}_{1}^{2p}-c^{3}\mathbf{D}_{1}^{3p}+\cdots+(-1)^{k}c^{k}\mathbf{D}_{1}^{kp}
\end{eqnarray*}
 with the largest integer $k$ that $kp\leq N$\@.
\end{corollary}
The following result, makes some savings in storage especially for
much larger matrices\@.
\begin{proposition}
Let $p\geq1$ and $0\leq i,j\leq N$, then 
\begin{eqnarray*}
 &  & (\mathbf{D}_{p})_{i,j}=(-1)^{p}(\mathbf{D}_{p})_{N-i,N-j},\,i,j\geq0.
\end{eqnarray*}
\end{proposition}
\begin{svmultproof}
For the nonzero elements of $\mathbf{D}_{p}$, from (\ref{eq3.4-1-1}),
we have
\begin{eqnarray*}
(\mathbf{D}_{p})_{N-i,N-j} & = & c_{p,N}\omega_{N-i,N-j}=c_{p,N}\sum_{k=\max(0,j-i)}^{\min(p,j-i+p)}{(-1)^{k}\frac{\binom{p}{k}\binom{N-p}{N-i-k}\binom{p}{(N-j)-(N-i)+k}}{\binom{N}{N-j}}}.
\end{eqnarray*}
To make the manipulations easier, we remove the max and min on the
summation limits (the only difference is that some zeros are added),
so
\begin{eqnarray*}
(\mathbf{D}_{p})_{N-i,N-j} & = & c_{p,N}\sum_{k=j-i}^{j-i+p}{(-1)^{k}\frac{\binom{p}{k}\binom{N-p}{N-i-k}\binom{p}{i-j+k}}{\binom{N}{j}}}\\
 & = & c_{p,N}\sum_{k=i-j}^{i-j+p}{(-1)^{k}\frac{\binom{p}{k-2i+2j}\binom{N-p}{N-k+i-2j}\binom{p}{k-i+j)}}{\binom{N}{j}}}.
\end{eqnarray*}
 Reversing the order of the last summation using $\sum_{k=m}^{n}a_{k}$=$\sum_{k=m}^{n}a_{m+n-k}$,
we get
\begin{eqnarray*}
(\mathbf{D}_{p})_{N-i,N-j} & = & c_{p,N}\sum_{k=i-j}^{i-j+p}{(-1)^{2i-2j+p-k}\frac{\binom{p}{p-k}\binom{N-p}{N-i-p+k}\binom{p}{i-j+p-k}}{\binom{N}{j}}}\\
 & = & c_{p,N}\sum_{k=i-j}^{i-j+p}{(-1)^{p-k}\frac{\binom{p}{k}\binom{N-p}{i-k}\binom{p}{j-i+k}}{\binom{N}{j}}}\\
 & = & (-1)^{p}c_{p,N}\sum_{k=\max(0,i-j)}^{\min(p,i-j+p)}{(-1)^{k}\frac{\binom{p}{k}\binom{N-p}{i-k}\binom{p}{j-i+k}}{\binom{N}{j}}}\\
 & = & (-1)^{p}c_{p,N}\sum_{k=\max(0,i-j)}^{\min(p,i-j+p)}{\mu_{i,k,j}}\\
 & = & (-1)^{p}c_{p,N}\omega_{i,j}=(-1)^{p}(\mathbf{D}_{p})_{i,j}.
\end{eqnarray*}
This completes the proof\@.
\end{svmultproof}

\begin{remark}
Corresponding to the non-orthogonal Bernstein  basis, the associated
orthonormalized basis $\{Q_{i,N}(x):\,i=0,\dots,N\}$ is obtained
by using the Gram-Schmidt algorithm\@. This basis fails to have properties
like (\ref{eq:firstDerCoef}) and (\ref{eq:SecDerCoeff})\@. To see
this, let $i=4$ and $N=4$\@. Then, it can be verified that 
\[
Q_{i,N}^{\prime}(x)=\sum_{i=0}^{N}c_{i}Q_{i,N}(x),
\]
in which all the coefficients $c_{0},\dots,c_{4}$ are nonzero\@.
It is also worth noting that the Chebyshev and Legendre polynomials
do not have properties like (\ref{eq:firstDerCoef}) and (\ref{eq:SecDerCoeff})\@.
In fact, the differentiation matrix of Chebyshev and also Legendre
polynomial basis of degree $N$ have ${\it [\frac{N}{2}]}$ nonzero
diagonals (see \citep{SaadatLegendre} and {[}\citealp{Shen},Sections
3.3 and 3.4{]} for more details)\@.
\end{remark}

\section{\label{sec:MainFormulationAndApplication}A new numerical method
for the time-fractional advection-dispersion equation}

This section is devoted to providing a numerical method using the
new Bernstein polynomial operational matrices with the Petrov-Galerkin
method for the problem (\ref{eq:main}) subject to the initial condition
given by (\ref{IV}) and the homogeneous boundary conditions (\ref{BVs})\@.

We rewrite the advection-dispersion equation (\ref{eq:main}) at $t=t_{k+1},\,k\geq0$
as
\begin{eqnarray*}
 &  & {D_{t}^{\alpha}}u(x,t_{k+1})=\kappa_{1}\frac{\partial^{2}}{\partial x^{2}}u(x,t_{k+1})-\kappa_{2}\frac{\partial}{\partial x}u(x,t_{k+1})+S(x,t_{k+1}),\quad x\in\Omega.
\end{eqnarray*}
Without loss of generality, we consider the problem (\ref{eq:main})-(\ref{BVs})
on the bounded domain $\Omega_{T}=[0,1]\times[0,T]$. Let $t_{k}=k\tau$,
$k=0,1,\dots,M$ where $\tau=\frac{T}{M}$ is the time step length.
Then, the time fractional derivative at time $t_{k+1}$ is approximated
as
\begin{eqnarray*}
 &  & {D_{t}^{\alpha}}u(x,t_{k+1})=\frac{1}{\Gamma(1-\alpha)}\sum_{j=0}^{k}\frac{u(x,t_{j+1})-u(x,t_{j})}{\tau}\int_{t_{j}}^{t_{j+1}}\frac{ds}{(t_{k+1}-s)^{\alpha}}+r_{\tau}^{k+1},\quad k\geq0,
\end{eqnarray*}
with the error term $r_{\tau}^{k+1}$\@. Consequently,
\begin{eqnarray}
 &  & {D_{t}^{\alpha}}u(x,t_{k+1})=\mu_{\tau}^{\alpha}\sum_{j=0}^{k}a_{k,j}^{\alpha}\left(u\left(x,t_{j+1}\right)-u\left(x,t_{j}\right)\right)+r_{\tau}^{k+1},\label{eq:4.1}
\end{eqnarray}
where $\mu_{\tau}^{\alpha}=\frac{1}{\tau^{\alpha}\Gamma(2-\alpha)}$
and $a_{kj}^{\alpha}=\left(k+1-j\right)^{1-\alpha}-\left(k-j\right)^{1-\alpha}$
for $k\geq0$ and $j=0,1,\dots,k$\@. A bound for the error is given
by
\begin{eqnarray}
|r_{\tau}^{k+1}| & \leq & \tilde{c}_{u}\tau^{2-\alpha},\label{eq:TempOrder2-alpha}
\end{eqnarray}
where the coefficient $\tilde{c}_{u}$ depends only on $u$ \citep{Deng}\@.
The scheme described here for time discretization is known as $L_{1}$
approximation, a common way in the numerical simulation of partial
fractional differential equations, for instance, see, \citep{Deng,Gao,Rame}\@.
Similar to \citep{Deng,Gao}, we define the time discrete fractional
differential operator $L_{t}^{\alpha}$ as
\begin{eqnarray}
 &  & L_{t}^{\alpha}u(x,t_{k+1})=\mu_{\tau}^{\alpha}\sum_{j=0}^{k}a_{k,j}^{\alpha}\left(u^{j+1}\left(x\right)-u^{j}\left(x\right)\right),\label{eq:Time}
\end{eqnarray}
in which $u^{j}(x)$ is an approximation for $u(x,t_{j})$\@. So
from (\ref{eq:4.1}), we have
\begin{eqnarray*}
 &  & {D_{t}^{\alpha}}u(x,t_{k+1})=L_{t}^{\alpha}u(x,t_{k+1})+r_{\tau}^{k+1}.
\end{eqnarray*}
Using $L_{t}^{\alpha}u(x,t_{k+1})$ as an approximation for ${D_{t}^{\alpha}}u(x,t_{k+1})$
in (\ref{eq:main}) yields the following time-discrete scheme:
\begin{eqnarray}
 &  & \mu_{\tau}^{\alpha}\sum_{j=0}^{k}a_{k,j}^{\alpha}\left(u^{j+1}\left(x\right)-u^{j}\left(x\right)\right)=\kappa_{1}\frac{\partial^{2}u^{k+1}(x)}{\partial x^{2}}-\kappa_{2}\frac{\partial u^{k+1}(x)}{\partial x}+S^{k+1}(x),\label{eq1}
\end{eqnarray}
for $k=0,...,M-1,$ from which we have
\begin{eqnarray}
 &  & \boldsymbol{L}u^{k+1}:=\mu_{\tau}^{\alpha}u^{k+1}(x)-\kappa_{1}\frac{\partial^{2}u^{k+1}(x)}{\partial x^{2}}+\kappa_{2}\frac{\partial u^{k+1}(x)}{\partial x}=f^{k+1}(x),\label{VIDE-1}
\end{eqnarray}
where $f^{k+1}(x)=\mu_{\tau}^{\alpha}\left(u^{k}(x)-\sum_{j=0}^{k-1}a_{k,j}^{\alpha}\left(u^{j+1}\left(x\right)-u^{j}\left(x\right)\right)\right)+S^{k+1}(x)$\@.
Note that $u^{0}(x)=g(x)$ is given by the initial condition (\ref{IV})
and the $f^{k+1}(x)$ is a known function at time step $t_{k+1}$\@.
For the error analysis , it will be useful to consider the homogeneous
case $S=0$, multiply the semidiscrete problem (\ref{VIDE-1}) by
$\tau^{\alpha}\Gamma(2-\alpha)$, rearrange the involving summation
and dropping $x$, to get \citep{Lin,Rame}
\begin{eqnarray}
 &  & u^{k+1}-\alpha_{1}\frac{\partial^{2}u^{k+1}}{\partial x^{2}}+\alpha_{2}\frac{\partial u^{k+1}}{\partial x}=(1-b_{1})u^{k}+\sum_{j=1}^{k-1}(b_{j}-b_{j+1})u^{k-j}+b_{k}u^{0},\quad k\geq1,\label{eq:SemiDiscrete}
\end{eqnarray}
and for $k=0$ 
\begin{eqnarray}
 &  & u^{1}-\alpha_{1}\frac{\partial^{2}u^{1}}{\partial x^{2}}+\alpha_{2}\frac{\partial u^{1}}{\partial x}=u^{0},\label{eq:SemiDescFork=00003D0}
\end{eqnarray}
where $\alpha_{i}=\kappa_{i}\tau^{\alpha}\Gamma(2-\alpha),\,i=1,2$
and 
\begin{eqnarray}
b_{j} & = & (j+1)^{1-\alpha}-j^{1-\alpha},\quad j=0,1,\dots,k.\label{eq:bj}
\end{eqnarray}
The boundary conditions are $u^{k+1}(0)=u^{k+1}(1)=0$ and the initial
condition is $u^{0}=g(x)$. By (\ref{eq:TempOrder2-alpha}), the error
for (\ref{eq:SemiDiscrete}), $r^{k+1}=\tau^{\alpha}\Gamma(2-\alpha)r_{\tau}^{k+1}$
is bounded as 
\begin{eqnarray}
 &  & |r^{k+1}|\leq c_{u}\tau^{2}.\label{eq:TruncationErrorAfterMultip}
\end{eqnarray}

To enforce the boundary conditions on the approximate solution, we
set $\mathbb{P}_{N}^{0}=\{\phi\in\mathbb{P}_{n}:\,\phi(0)=\phi(1)=0\}$
with $dim(\mathbb{P}_{N}^{0})=N-1$\@. Define $\tilde{\Phi}=[\phi_{i}(x)=B_{i,N}(x):\,i=1,\dots,N-1]^{T},$
where from now on, we set $a=0$ and $b=1$\@. Since $\tilde{\Phi}$
is a basis for $X_{N}$, we can expand the approximate solution at
time step $t=t_{k+1}$ as
\begin{eqnarray}
 &  & u^{k+1}(x)\approx u_{N}^{k+1}(x)=\sum_{i=1}^{N-1}{c_{i}^{k+1}B_{i,N}(x)}=(\boldsymbol{c}^{k+1})^{T}\tilde{\phi},\quad k\geq0,\label{eq:Expansion}
\end{eqnarray}
 where $\boldsymbol{c}^{k+1}=[c_{i}^{k+1}:\,i=1,\dots,N-1]^{T}$\@.
It is easy to see that the associated differentiation matrices are
given as
\begin{eqnarray}
 &  & \tilde{\Phi}^{\prime}=\tilde{\mathbf{D}}_{1}\tilde{\Phi}+\boldsymbol{d}_{1},\label{eq:MatOper1}\\
 &  & \tilde{\Phi}^{\prime\prime}=\tilde{\mathbf{D}}_{2}\tilde{\Phi}+\boldsymbol{d}_{2},\label{eq:Mat.Oper2}
\end{eqnarray}
where the tridiagonal and pentadiagonal matrices $\tilde{\mathbf{D}}_{1}$
and $\tilde{\mathbf{D}}_{2}$ are extracted respectively by removing
the first and last rows and the first and last columns of $\mathbf{D}_{1}$
and $\mathbf{D}_{2}$\@. The $(N-1)$-dimensional vectors $\boldsymbol{d}_{1}$
and $\boldsymbol{d}_{2}$ are given by $\boldsymbol{d}_{1}=N[\phi_{0},0,\ldots.,0,-\phi_{N}]^{T}$
and $\boldsymbol{d}_{2}=N(N-1)[-2\phi_{0},\phi_{0},0,\ldots.,0,\phi_{N}-2\phi_{N}]^{T}$\@.
By substituting (\ref{eq:Expansion}) for $u^{k+1}$ in (\ref{VIDE-1})
and using (\ref{eq:MatOper1}) and (\ref{eq:Mat.Oper2}), we get
\begin{eqnarray*}
 &  & \boldsymbol{L}u_{N}^{k+1}(x)=(\boldsymbol{c}^{k+1})^{T}(\mu_{\tau}^{\alpha}\tilde{\Phi}-\kappa_{1}(\tilde{\mathbf{D}}_{2}\tilde{\Phi}+\boldsymbol{d}_{2})+\kappa_{2}(\tilde{\mathbf{D}}_{1}\tilde{\Phi}+\boldsymbol{d}_{1})),
\end{eqnarray*}
and the residual is given by 
\begin{eqnarray*}
 &  & \boldsymbol{R}_{N}^{k+1}(x):=\boldsymbol{L}u_{N}^{k+1}(x)-f^{k+1}(x).
\end{eqnarray*}
To obtain the expansion coefficients $\boldsymbol{c}_{k+1}$, we enforce
the residual to be orthogonal to the dual basis functions $\psi_{j}:=B_{j,N}^{\star},\,j=1,\dots,N-1$
given by (\ref{eq:Dual}), i.e., 
\begin{eqnarray}
 &  & ((\boldsymbol{c}^{k+1})^{T}(\mu_{\tau}^{\alpha}\tilde{\Phi}-\kappa_{1}(\tilde{\mathbf{D}}_{2}\tilde{\Phi}+\boldsymbol{d}_{2})+\kappa_{2}(\tilde{\mathbf{D}}_{1}\tilde{\Phi}+\boldsymbol{d}_{1}))-f^{k+1},\,\psi_{j})=0,\quad j=1,\dots,N-1,\label{eq:LinearSys}
\end{eqnarray}
where $(\cdot,\cdot)$ is the standard $L^{2}-$inner product\@.
Define the dual basis vector $\tilde{\Psi}^{T}=[\psi_{i}(x):\,i=1,\dots,N-1]$\@.
Using the biorthogonal identity (\ref{eq2.11}), we have $\tilde{\Phi}\Psi^{T}=\mathbf{I}_{N-1},$
the identity matrix of order $N-1$\@. Since $(\boldsymbol{d}_{1},\psi_{j})=0$
and $(\boldsymbol{d}_{2},\psi_{j})=0$ for $0<j<N,$ the above linear
system can be written as a pentadiagonal linear system as follows:
\begin{eqnarray}
 &  & \mathbf{A}\boldsymbol{c}^{k+1}=\boldsymbol{b}^{k+1},\quad k=0,1,\dots,M-1,\label{eq:mainMatrixSystem}
\end{eqnarray}
where
\begin{eqnarray}
 &  & \mathbf{A}^{T}=\mu_{\tau}^{\alpha}\mathbf{I}_{N-1}-\kappa_{1}\tilde{\mathbf{D}}_{2}+\kappa_{2}\tilde{\mathbf{D}}_{1},\label{eq:matrixAtrans}
\end{eqnarray}
 and $\boldsymbol{b}^{k+1}=[(f^{k+1},\psi_{i}):\,i=1,\dots,N-1]^{T}$\@.
The integrals in $\boldsymbol{b}^{k+1}$ may be computed by a numerical
method like Gauss quadrature rules\@. In our computational test examples,
we use 20 point Gauss-Legendre rule on the unit interval\@.

It is worth to note that the coefficient matrix of the linear system
(\ref{eq:mainMatrixSystem}) is independent of the time step $k$.
So for a fixed $N$ and $M$, its $LU$-decomposition is performed
just once and is used for all time steps\@. The LU-decomposition
for a banded matrix with $2p+1$-bandwidth is done just by $O(Np^{2}$)
arithmetic operations and the number of operations to modify the right
hand side and performing back substitution is $O(Np)$ {[}\citealp{Golub},
Section 4.3{]}. So obtaining the numerical solution of the problem
(\ref{eq:main}) on the bounded domain $[0,1]\times[0,T]$ with the
proposed method requires just $O(MNp(p+1))$ arithmetic operations.

Starting from $u_{N}^{0}=g(x)$, the approximate solution (\ref{eq:Expansion})
at $t_{k+1},\,k\geq0$ is obtained by solving the linear system (\ref{eq:mainMatrixSystem}). 

The following result states the condition required for the existence
of the solution for the linear system (\ref{eq:mainMatrixSystem}).
\begin{lemma}
Let $\left\Vert \cdot\right\Vert $ be an induced matrix norm and
suppose that the time step length $\tau$ is such that
\begin{eqnarray}
 &  & \kappa_{1}\|\tilde{\mathbf{D}}_{2}\|+\kappa_{2}\|\tilde{\mathbf{D}}_{1}\|<\mu_{\tau}^{\alpha}.\label{eq:Condition}
\end{eqnarray}
 Then, the coefficient matrix $\mathbf{A}$ is nonsingular, i.e.,
the linear system (\ref{eq:mainMatrixSystem}) has a unique solution.
\end{lemma}
\begin{svmultproof}
Assume that (\ref{eq:Condition}) holds and $\mathbf{A}$ is singular.
Then, there is a nonzero vector $x$ with $\|x\|=1$ such that $\mathbf{A}^{T}x=0$
and so $\mu_{\tau}^{\alpha}=\|(\kappa_{1}\tilde{\mathbf{D}}_{2}-\kappa_{2}\tilde{\mathbf{D}}_{1})x\|.$
Then, $\mu_{\tau}^{\alpha}\leq\kappa_{1}\|\tilde{\mathbf{D}}_{2}\|+\kappa_{2}\|\tilde{\mathbf{D}}_{1}\|$
which contradicts (\ref{eq:Condition}).
\end{svmultproof}

Note that $\mu_{\tau}^{\alpha}\rightarrow\infty$ as time mesh size
tends to zero. So the principal significance of the lemma is that
it allows one to choose time mesh size $\tau$ small enough in order
to satisfy the condition (\ref{eq:Condition}). However, by Maple,
we found that $\det\mathbf{A}\neq0$ for $1\leq N\leq100$ without
imposing condition (\ref{eq:Condition}) on the parameters, except
for advection and dispersion coefficients $\kappa_{1}>0$ and $\kappa_{2}>0$
which are positive due to the nature of the problem, and it can be
written as 
\begin{eqnarray*}
 &  & \det{\mathbf{A}}=\sum_{i=0}^{N-1}{\sum_{j=0}^{[\frac{N-i-1}{2}]}{d_{i,j}\kappa_{2}^{2j}}\kappa_{1}^{i}}>0,
\end{eqnarray*}
with all positive coefficients $d_{i,j}$.

It is easy to verify that $\|\tilde{\mathbf{D}}_{1}\|_{\infty}=2(N-1)$
for $N\geq4$ and
\begin{eqnarray*}
 &  & \|\tilde{\mathbf{D}}_{2}\|_{\infty}=\begin{cases}
\begin{array}{cc}
2N^{2}+10N-68, & 7\leq N\leq15,\\
4N^{2}-28N+40, & N\geq16.
\end{array}\end{cases}
\end{eqnarray*}
and therefore, from (\ref{eq:matrixAtrans}), we have 
\begin{eqnarray*}
 &  & \|\mathbf{A\|_{\infty}}\leq\begin{cases}
\begin{array}{cc}
\mu+\kappa_{1}(2N^{2}+10N-68)+2\kappa_{2}(N-1), & 7\leq N\leq15,\\
\mu+\kappa_{1}(4N^{2}-28N+40)+2\kappa_{2}(N-1), & N\geq16,
\end{array}\end{cases}
\end{eqnarray*}
but since we do not have a compact formula for $\|\mathbf{A}^{-1}\|_{\infty}$,
we only provide the condition number of the coefficient matrix $\mathbf{A}$
for some numeric values.

In Table \ref{tab:Cond}, we provide the condition number of $\mathbf{A}$
and compare it with the condition number of the Hilbert matrix $\mathbf{H}$,
with respect to the infinity norm, i.e., $C_{\infty}(\mathbf{A}):=Cond_{\infty}(\mathbf{A})=\|\mathbf{A}\|_{\infty}\|\mathbf{A}^{-1}\|_{\infty}$
for $\alpha=0.5,\,\tau=1/40$ for numeric values of $\kappa_{1}$
and $\kappa_{2}$. It shows that the condition number of $\mathbf{A}$
is relatively small and the ratio $R_{\infty}=C_{\infty}(\mathbf{A})/C_{\infty}(\mathbf{H})\rightarrow0$
as $N\rightarrow\infty$ where $N$ is the size of the matrix. So
we can expect good numerical results. This is supported by the numerical
examples in the next section. On the other hand, as we mentioned before
the transformation between Bernstein and Power basis is ill-conditioned
and should be avoided for large $N$. In \citep{Farouki1999}, Farouki
showed that the transformation matrix has condition number $C_{\infty}\approx3^{N+1}\sqrt{(N+1)/4\pi}$.
\begin{flushleft}
\begin{table}[h]
\caption{\label{tab:Cond}Condition number of $\mathbf{A}$ and comparison
with the condition number of the Hilbert matrix}
\noindent\resizebox{\textwidth}{!}{
\begin{tabular}{c|ccccccccc}
\hline 
\multicolumn{1}{c}{} & $N$ & 4 & 5 & 6 & 7 & 8 & 9 & 10 & 11\tabularnewline
\hline 
\multirow{2}{*}{$\begin{array}{c}
\kappa_{1}=0.1,\\
\kappa_{2}=2
\end{array}$} & $C_{\infty}(\mathbf{A})$ & 5.31 & 8.03 & 12.90 & 27.41 & 54.77 & 100.74 & 210.08 & 463.47\tabularnewline
 & $R_{\infty}$ & 1.87E-04 & 8.51E-06 & 4.44E-07 & 2.78E-08 & 1.62E-09 & 9.16E-11 & 5.94E-12 & 3.76E-13\tabularnewline
\hline 
\multirow{2}{*}{$\begin{array}{c}
\kappa_{1}=1,\\
\kappa_{2}=1
\end{array}$} & $C_{\infty}(\mathbf{A})$ & 1.57 & 1.73 & 1.86 & 7.33 & 11.76 & 19.78 & 34.82 & 63.57\tabularnewline
 & $R_{\infty}$ & 5.54E-05 & 1.84E-06 & 6.43E-08 & 7.45E-09 & 3.47E-10 & 1.80E-11 & 9.85E-13 & 5.15E-14\tabularnewline
\hline 
\end{tabular}}
\end{table}
\par\end{flushleft}

\section{\label{sec:Error}Error estimation}

\subsection{Stability and convergence of the semidiscrete scheme}

We will carry out the error estimation for the homogeneous case of
(\ref{eq:main}), i.e., for $S=0$. 

Due to the presence of the first-order advection term, it is convenient
for the error analysis to multiply both sides of (\ref{eq:SemiDiscrete})
by an integrating factor and use a weighted variational formulation
(see e.g., {[}\citealp{Shen}; Section 4.4{]}). However, in order
to utilize the biorthogonality (\ref{eq2.11}) providing banded sparse
linear systems the matrix, the formulation of our method in Section
\ref{sec:MainFormulationAndApplication} and the numerical computations
in Section \ref{sec:Numerical-examples}) are presented with the spectral
formulation without weight function.

Multiplying (\ref{eq:SemiDiscrete}) by $w=\exp(-\frac{\alpha_{2}}{\alpha_{1}}x)=\exp(-\frac{\kappa_{2}}{\kappa_{1}}x)$,
the equation (\ref{eq:SemiDiscrete}) is written as 
\begin{eqnarray}
 &  & wu^{k+1}-\alpha_{1}\frac{\partial}{\partial x}\left(w\frac{\partial u^{k+1}}{\partial x}\right)=w(1-b_{1})u^{k}+\sum_{j=1}^{k-1}(b_{j}-b_{j+1})wu^{k-j}+b_{k}wu^{0}.\label{eq:SemiDiscreteNew}
\end{eqnarray}
The variational formulation is then written as
\begin{eqnarray}
 &  & (u^{k+1},v)_{w}+\alpha_{1}(\frac{\partial u^{k+1}}{\partial x},\frac{\partial v}{\partial x})_{w}=(1-b_{1})(u^{k},v)_{w}\nonumber \\
 &  & \quad\quad\quad+\sum_{j=1}^{k-1}(b_{j}-b_{j+1})(u^{k-j},v)_{w}+b_{k}(u^{0},v)_{w},\quad k\geq1,\,\forall v\in H_{0,w}^{1}(\Omega),\label{eq:weakFormulationNew}
\end{eqnarray}
and for $k=0,$
\begin{eqnarray}
 &  & (u^{1},v)_{w}+\alpha_{1}(\frac{\partial u^{1}}{\partial x},\frac{\partial v}{\partial x})_{w}=(u^{0},v)_{w},\quad\forall v\in H_{0,w}^{1}(\Omega).\label{eq:weakFormulationNew,k=00003D0}
\end{eqnarray}
Accordingly the Galerkin spectral discretization is to find $u_{N}^{k+1}\in\mathbb{P}_{N}^{0}(\Omega)=H_{0,w}^{1}(\Omega)\cap\mathbb{P}_{N}(\Omega)$
such that
\begin{eqnarray}
 &  & (u_{N}^{k+1},v_{N})_{w}+\alpha_{1}(\frac{\partial u_{N}^{k+1}}{\partial x},\frac{\partial v_{N}}{\partial x})_{w}=(1-b_{1})(u_{N}^{k},v_{N})_{w}\nonumber \\
 &  & \quad\quad\quad+\sum_{j=1}^{k-1}(b_{j}-b_{j+1})(u_{N}^{k-j},v_{N})_{w}+b_{k}(u_{N}^{0},v_{N})_{w},\quad k\geq1,\,\forall v_{N}\in\mathbb{P}_{N}^{0}(\Omega).\label{eq:SpectralDiscretNew}
\end{eqnarray}
for $k=0,$
\begin{eqnarray}
 &  & (u_{N}^{1},v_{N})_{w}+\alpha_{1}(\frac{\partial u_{N}^{1}}{\partial x},\frac{\partial v_{N}}{\partial x})_{w}=(u_{N}^{0},v_{N})_{w},\quad\forall v_{N}\in H_{0,w}^{1}(\Omega).\label{eq:SpectralDiscretNew,k=00003D0}
\end{eqnarray}
We define the following inner product and the associated energy norm
on $H_{0,w}^{1}(\Omega)$: 
\begin{eqnarray}
 &  & (u,v)_{w}=\int_{\Omega}uvwd\Omega,\quad(u,v)_{1,w}=(u,v)_{w}+\alpha_{1}(\frac{\partial u}{\partial x},\frac{\partial v}{\partial x})_{w},\quad\|u\|_{1,w}=(u,u)_{1,w}^{\frac{1}{2}}.\label{eq:EnergyNormWeighted}
\end{eqnarray}
It is worth noting that from $\kappa_{2}\leq\kappa_{1}$ as it happens
for real advection diffusion problems, we get $0\leq w(x)\leq1$ for
$x\in\Omega$, so $L^{2}(\Omega)\subseteq L_{w}^{2}(\Omega)$ and
$H^{1}(\Omega)\subseteq H_{w}^{1}(\Omega)$. 

The following result presents the unconditional stability of the the
scheme (\ref{eq:weakFormulationNew}). 
\begin{theorem}
\label{Th:stabilityForThe-weak-semidiscrete}The weak semidiscrete
scheme (\ref{eq:weakFormulationNew}) is unconditionally stable: 
\begin{eqnarray}
 &  & \|u^{k}\|_{1,w}\leq\|u^{0}\|_{w},\quad k=1,\dots,M.\label{eq:stabilityOfSemiDiscrete}
\end{eqnarray}
\end{theorem}
\begin{svmultproof}
Let $v=u^{1}$ in (\ref{eq:weakFormulationNew,k=00003D0}). Then,
\begin{eqnarray*}
 &  & (u^{1},u^{1})_{w}+\alpha_{1}(\frac{\partial u^{1}}{\partial x},\frac{\partial u^{1}}{\partial x})_{w}=(u^{0},u^{1})_{w}.
\end{eqnarray*}
Using the Schwarz inequality, the inequality $\|v\|_{w}\leq\|v\|_{1,w}$,
and dividing both sides by $\|u^{1}\|{}_{1,w}$, one immediately gets
(\ref{eq:stabilityOfSemiDiscrete}) for $k=1$. Now suppose that (\ref{eq:stabilityOfSemiDiscrete})
holds for $k\leq n.$ Taking $v=u^{n+1}$ in (\ref{eq:weakFormulationNew}),
we get
\begin{eqnarray*}
 &  & (u^{n+1},u^{n+1})_{w}+\alpha_{1}(\frac{\partial u^{n+1}}{\partial x},\frac{\partial u^{n+1}}{\partial x})_{w}=(1-b_{1})(u^{n},u^{n+1})_{w}\\
 &  & \quad\quad\quad\hfill+\sum_{j=1}^{n-1}(b_{j}-b_{j+1})(u^{n-j},u^{n+1})_{w}+b_{n}(u^{0},u^{n+1})_{w}.
\end{eqnarray*}
Note that the sequence $\{b_{j}\}$ defined in (\ref{eq:bj}) is decreasing
and converging to zero with $b_{0}=1$. Hence the RHS coefficients
in equation (\ref{eq:SemiDiscrete}) are positive. Again using $\|v\|_{w}\leq\|v\|_{1,w}$,
the Schwarz inequality and dividing both sides by $\|u^{n+1}\|_{1,w}$,
we get 
\begin{eqnarray*}
\|u^{n+1}\|_{1,w} & \leq & (1-b_{1})\|u^{n}\|_{w}+\sum_{j=1}^{n-1}(b_{j}-b_{j+1})\|u^{n-j}\|_{w}+b_{n}\|u^{0}\|_{w}\\
 & \leq & \left((1-b_{1})+\sum_{j=1}^{n-1}(b_{j}-b_{j+1})+b_{n}\right)\|u^{0}\|_{w}=\|u^{0}\|_{w},
\end{eqnarray*}
that is (\ref{eq:stabilityOfSemiDiscrete}) for $k=n+1$. This completes
the proof. 
\end{svmultproof}

\begin{theorem}
\label{Th:ConvergenceForSemidiscrete}Let $u(x,t)$ be the exact solution
of the problem (\ref{eq:main}) with the initial and boundary conditions
(\ref{IV})-(\ref{BVs}) and $u^{k},\,k=1,\dots,M,$ be the solution
of the the semidiscrete problem (\ref{eq:SemiDiscrete}). Then,
\begin{eqnarray}
 &  & \|u(t_{k})-u^{k}\|_{1,w}\leq\frac{c_{u}}{1-\alpha}T^{\alpha}\tau^{2-\alpha},\quad0<\alpha<1,\label{eq:Error u(t_k)-u^k}\\
 &  & \|u(t_{k})-u^{k}\|_{1,w}\leq c_{u}T\tau,\quad\qquad\quad\,as\,\alpha\rightarrow1.\label{eq:Error u-u^k alpha->1}
\end{eqnarray}
\end{theorem}
\begin{svmultproof}
We first prove that 
\begin{eqnarray}
 &  & \|u(t_{k})-u^{k}\|_{1,w}\leq\frac{c_{u}}{b_{k-1}}\tau^{2},\quad k=1,\dots,M.\label{eq:Auxilary}
\end{eqnarray}
Define $e^{k}:=u(t_{k})-u^{k}.$ By (\ref{eq:main}) and (\ref{eq:SemiDescFork=00003D0}),
we derive
\begin{eqnarray*}
 &  & (e^{1},v)_{w}+\alpha_{1}(\frac{\partial e^{1}}{\partial x},\frac{\partial v}{\partial x})_{w}=(e^{0},v)_{w}+(r^{1},v)_{w},\quad\forall v\in H_{0,w}^{1}(\Omega).
\end{eqnarray*}
Let $v=e^{1}$, then using $e^{0}=0$, $\|v\|_{w}\leq\|v\|_{1,w}$
and (\ref{eq:TruncationErrorAfterMultip}), we obtain
\begin{eqnarray}
 &  & \|e^{1}\|_{1,w}\leq c_{u}\tau^{2},\label{eq:Forj=00003D1}
\end{eqnarray}
that is (\ref{eq:Auxilary}) for $k=1$. Now suppose that (\ref{eq:Auxilary})
holds for $k=1,\dots,n$. Using (\ref{eq:main}) and (\ref{eq:SemiDiscrete}),
we get
\begin{eqnarray*}
 &  & (e^{n+1},v)_{w}+\alpha_{1}(\frac{\partial e^{n+1}}{\partial x},\frac{\partial v}{\partial x})_{w}=(1-b_{1})(e^{n},v)_{w}\\
 &  & \quad\quad\quad+\sum_{j=1}^{n-1}(b_{j}-b_{j+1})(e^{n-j},v)_{w}+b_{n}(e^{0},v)_{w}+(r^{n+1},v)_{w},\quad\forall v\in H_{0,w}^{1}(\Omega).
\end{eqnarray*}
Taking $v=e^{n+1}$, we have
\begin{eqnarray*}
\|e^{n+1}\|_{1,w}^{2} & \leq & (1-b_{1})\|e^{n}\|_{w}\|e^{n+1}\|_{1,w}+\sum_{j=1}^{n-1}(b_{j}-b_{j+1})\|e^{n-j}\|_{w}\|e^{n+1}\|_{1,w}+\|r^{n+1}\|_{w}\|e^{n+1}\|_{1,w},\\
\Rightarrow\|e^{n+1}\|_{1,w} & \leq & (1-b_{1})\frac{c_{u}}{b_{n-1}}\tau^{2}+\sum_{j=1}^{n-1}(b_{j}-b_{j+1})\frac{c_{u}}{b_{n-j-1}}\tau^{2}+c_{u}\tau^{2}\\
 & \leq & \left((1-b_{1})+\sum_{j=1}^{n-1}(b_{j}-b_{j+1})+b_{n}\right)\frac{c_{u}}{b_{n}}\tau^{2}=\frac{c_{u}}{b_{n}}\tau^{2},
\end{eqnarray*}
i.e., (\ref{eq:Auxilary}) holds for $k=n+1$. This proves the (\ref{eq:Auxilary}). 

Applying the mean value theorem to the function $f(t)=t^{1-\alpha}$,
$\exists\xi$ $k-1<\xi<k\leq M$ such that
\begin{eqnarray*}
 &  & b_{k-1}\tau^{-\alpha}=\frac{(k\tau)^{1-\alpha}-(\tau(k-1))^{1-\alpha}}{\tau}=(1-\alpha)(\xi\tau)^{-\alpha}\geq(1-\alpha)(k\tau)^{-\alpha}\geq(1-\alpha)(T)^{-\alpha},
\end{eqnarray*}
so we obtain 
\begin{eqnarray}
 &  & \frac{c_{u}}{b_{k-1}}\tau^{2}\leq\frac{c_{u}}{1-\alpha}T^{\alpha}\tau^{2-\alpha}.\label{eq:boundFor b_k}
\end{eqnarray}
This together with (\ref{eq:Auxilary}) gives (\ref{eq:Error u(t_k)-u^k}). 

To prove (\ref{eq:Error u-u^k alpha->1}), we first prove 
\begin{eqnarray}
 &  & \|u(t_{k})-u^{k}\|_{1,w}\leq c_{u}k\tau^{2},\quad k=1,\dots,M.\label{eq:Auxilary2}
\end{eqnarray}
 From (\ref{eq:Forj=00003D1}), the relation (\ref{eq:Auxilary2})
holds for $k=1$. Let (\ref{eq:Auxilary2}) holds for $k=1,\dots,n,\:n\leq M-1$.
Then, from (\ref{eq:main}), (\ref{eq:SemiDiscrete}) and (\ref{eq:TruncationErrorAfterMultip}),
we obtain
\begin{eqnarray*}
\|e^{n+1}\|_{1,w} & \leq & (1-b_{1})\|e^{n}\|_{w}+\sum_{j=1}^{n-1}(b_{j}-b_{j+1})\|e^{n-j}\|_{w}+\|r^{n+1}\|_{w}\\
 & \leq & \left((1-b_{1})\frac{n}{n+1}+\sum_{j=1}^{n-1}(b_{j}-b_{j+1})\frac{n-j}{n+1}+\frac{1}{(n+1)}\right)c_{u}(n+1)\tau^{2}\\
 & \leq & \left((1-b_{1})\frac{n}{n+1}+(b_{1}-b_{n})\frac{n}{n+1}-(b_{1}-b_{n})\frac{1}{n+1}+\frac{1}{(n+1)}\right)c_{u}(n+1)\tau^{2}\\
 & = & \left(1-b_{n}\frac{n}{n+1}-(b_{1}-b_{n})\frac{1}{n+1}\right)c_{u}(n+1)\tau^{2}\leq c_{u}(n+1)\tau^{2}.
\end{eqnarray*}
This gives (\ref{eq:Auxilary2}) for $k=n+1$. The proof for (\ref{eq:Auxilary2})
is done. Now the relation (\ref{eq:Auxilary2}) with $k\tau\leq T$
gives (\ref{eq:Error u-u^k alpha->1}).
\end{svmultproof}

\subsection{Convergence of the full discretization scheme}

Let $\pi_{N,w}^{1,0}$ be the $H^{1}$-orthogonal projection operator
from $H_{0,w}^{1}(\Omega)$ into $\mathbb{P}_{N}^{0}(\Omega)$ related
to the energy norm $\|\cdot\|_{1,w}$. We have the error estimation
by using {[}\citealp{Lin}; Relation (4.3){]}
\begin{eqnarray}
 &  & \|u-\pi_{N,w}^{1,0}u\|_{1,w}\leq cN^{1-m}\|u\|_{w,m},\quad u\in H_{0,w}^{m}(\Omega)\cap H_{0,w}^{1}(\Omega),\,m\geq1.\label{eq:projectionError}
\end{eqnarray}

The idea for the proof of the following result comes from the paper
\citep{Lin} in which the authors use a collocation spectral technique
for the subdiffusion equations.
\begin{theorem}
Let $\{u^{k}\}_{k=0}^{M}$ be the solution of (\ref{eq:weakFormulationNew})
and $\{u_{N}^{k}\}_{k=0}^{M}$ be the solution of the spectral discretization
(\ref{eq:SpectralDiscretNew}). Assume $u^{0}=\pi_{N,w}^{1,0}u^{0}$
and $u^{k}\in H_{w}^{m}(\Omega)\cap H_{0,w}^{1}(\Omega)$ for some
$m>1.$ Then,
\begin{eqnarray}
 &  & \|u^{k}-u_{N}^{k}\|_{1,w}\leq\frac{c}{1-\alpha}\tau^{-\alpha}N^{1-m}\max_{0\leq j\leq k}\|u^{j}\|_{m,w},\quad0<\alpha<1,\nonumber \\
 &  & \|u^{k}-u_{N}^{k}\|_{1,w}\leq cN^{1-m}\sum_{j=0}^{k}\|u^{j}\|_{m,w},\hfill\qquad\,\quad\alpha\rightarrow1,\label{eq:Error u^k-(u_N)^k}
\end{eqnarray}
for $k=1,\dots,M$, where $c$ depends only on $T^{\alpha}$.
\end{theorem}
\begin{svmultproof}
By definition (\ref{eq:EnergyNormWeighted}), we have from $(u^{k+1}-\pi_{N,w}^{1,0}u^{k+1},v_{N})_{1,w}=0,\,\forall v_{N}\in\mathbb{P}_{N}^{0}(\Omega)$
that 
\begin{eqnarray*}
 &  & (\pi_{N,w}^{1,0}u^{k+1},v_{N})_{w}+\alpha_{1}(\frac{\partial\pi_{N,w}^{1,0}u^{k+1})}{\partial x},\frac{\partial v_{N}}{\partial x})_{w}=(u^{k+1},v_{N})_{w}+\alpha_{1}(\frac{\partial u^{k+1}}{\partial x},\frac{\partial v_{N}}{\partial x})_{w},\quad\forall v_{N}\in\mathbb{P}_{N}^{0}(\Omega),
\end{eqnarray*}
Now using (\ref{eq:SemiDiscrete}), we get
\begin{eqnarray}
 &  & (\pi_{N,w}^{1,0}u^{k+1},v_{N})_{w}+\alpha_{1}(\frac{\partial\pi_{N,w}^{1,0}u^{k+1})}{\partial x},\frac{\partial v_{N}}{\partial x})_{w}=(1-b_{1})(u^{k},v_{N})_{w}\nonumber \\
 &  & \quad\quad\quad+\sum_{j=1}^{k-1}(b_{j}-b_{j+1})(u^{k-j},v_{N})_{w}+b_{k}(u^{0},v_{N})_{w},\quad\forall v_{N}\in\mathbb{P}_{N}^{0}(\Omega).\label{eq:11}
\end{eqnarray}
Defining the errors $e_{N}^{k+1}=u^{k+1}-u_{N}^{k+1}$ and $\tilde{e}_{N}^{k+1}=\pi_{N,w}^{1,0}u^{k+1}-u_{N}^{k+1}$
and subtracting (\ref{eq:11}) from (\ref{eq:SpectralDiscretNew}),
we have
\begin{eqnarray*}
 &  & (\tilde{e}_{N}^{k+1},v_{N})_{w}+\alpha_{1}(\frac{\partial\tilde{e}_{N}^{k+1}}{\partial x},\frac{\partial v_{N}}{\partial x})_{w}=(1-b_{1})(e_{N}^{k},v_{N})_{w}\\
 &  & \quad\quad\quad+\sum_{j=1}^{k-1}(b_{j}-b_{j+1})(e_{N}^{k-j},v_{N})_{w}+b_{k}(e_{N}^{0},v_{N})_{w},\quad\forall v_{N}\in\mathbb{P}_{N}^{0}(\Omega),
\end{eqnarray*}
Hence,
\begin{eqnarray*}
 &  & \|\tilde{e}_{N}^{k+1}\|_{1,w}\leq(1-b_{1})\|e_{N}^{k}\|_{w}+\sum_{j=1}^{k-1}(b_{j}-b_{j+1})\|e_{N}^{k-j}\|_{w}+b_{k}\|e_{N}^{0}\|_{w}.
\end{eqnarray*}
Now using $\|e_{N}^{k+1}\|_{1,w}\leq\|\tilde{e}_{N}^{k+1}\|_{1,w}+\|u^{k+1}-\pi_{N,w}^{1,0}u^{k+1}\|_{1,w}$,
we have
\begin{eqnarray*}
 &  & \|e_{N}^{k+1}\|_{1,w}\leq(1-b_{1})\|e_{N}^{k}\|_{w}+\sum_{j=1}^{k-1}(b_{j}-b_{j+1})\|e_{N}^{k-j}\|_{w}+b_{k}\|e_{N}^{0}\|_{w}+cN^{1-m}\|u^{k+1}\|_{m}.
\end{eqnarray*}
As in the proof of Theorem \ref{Th:ConvergenceForSemidiscrete}, it
is first proved by induction that:
\begin{eqnarray*}
 &  & \|e_{N}^{k+1}\|_{1,w}\leq\frac{1}{b_{k-1}}\max_{0\leq j\leq k}\|u^{j}-\pi_{N,w}^{1,0}u^{j}\|_{1,w},\quad0<\alpha<1,\\
 &  & \|e_{N}^{k+1}\|_{1,w}\leq\sum_{j=0}^{k}\|u^{j}-\pi_{N,w}^{1,0}u^{j}\|_{1,w},\quad\alpha\rightarrow1,
\end{eqnarray*}
 for $0\leq k\leq M.$ Then, using (\ref{eq:boundFor b_k}) and the
error bound (\ref{eq:projectionError}) the desired result is obtained. 
\end{svmultproof}

The following theorem is obtained by $||u(\cdot,t_{k})-u_{N}^{k}||_{1,w}\leq||u(\cdot,t_{k})-u^{k}||_{1,w}+||u^{k}-u_{N}^{k}||_{1,w}$
with (\ref{eq:Error u(t_k)-u^k}) and (\ref{eq:Error u^k-(u_N)^k})
for the first and second term of RHS, respectively.
\begin{theorem}
Let $u$ be the solution of equation (\ref{eq:main}) with conditions
(\ref{IV})-(\ref{BVs}) and $u_{N}^{k}$ be the solution of (\ref{eq:SpectralDiscretNew})
with $u_{N}^{0}=\pi_{N,w}^{1,0}u^{0}.$ If $u\in H_{w}^{m}(\Omega)\cap H_{0,w}^{1}(\Omega)$,
then
\begin{eqnarray}
 &  & \|u(t_{k})-u_{N}^{k}\|_{1,w}\leq\frac{CT^{\alpha}}{1-\alpha}(c_{u}\tau^{2-\alpha}+c\tau^{-\alpha}N^{1-m}\sup_{t\in(0,T)}\|u(x,t)\|_{m,w}),\quad k\leq M,\,0<\alpha<1,\label{eq:SpatialRate}\\
 &  & \|u(t_{k})-u_{N}^{k}\|_{1,w}\leq T^{\alpha}(c_{u}\tau+c\tau^{-1}N^{1-m}\sup_{t\in(0,T)}\|u(x,t)\|_{m,w}),\quad k\leq M,\quad\alpha\rightarrow1.\nonumber 
\end{eqnarray}
where $C$ and $c$ are constants independent of $N,\tau,T$ and $c_{u}$
depends only on $u$.
\end{theorem}
It is seen that the temporal and spatial rate of convergence are $O(\tau^{2-\alpha})$
and $O(N^{1-m})$, respectively, where $m$ is an index of regularity
of the underlying function. In the next section, Figures \ref{fig:The-spectral-convergence}
and \ref{figEx4} are provided to show the spectral accuracy of the
method in space and Table \ref{tab:TemporalOrderU=00003DxCos(..}
to show the rate in time for some numerical tests.

\section{\label{sec:Numerical-examples}Numerical examples}

In this section, we provide some numerical examples to show the efficiency
and accuracy of the method. We use the discrete $L^{2}$ and $L^{\infty}$
error measures at time $t_{M}=T=1$ as 
\begin{eqnarray}
 &  & L_{2}^{T}:=\left(\int_{0}^{1}{|u(x,T)-U(x,T)|^{2}dx}\right)^{1/2}\approx\left(\frac{1}{\mathcal{N}}\sum_{j=0}^{\mathcal{N}-1}{|u(x_{j},T)-U(x_{j},T)|^{2}}\right)^{1/2},\label{eq:L_2}\\
 &  & L_{\infty}^{T}:=\max_{0\leq x\leq1}{|u(x,T)-U(x,T)|}\approx\max_{0\leq j\leq\mathcal{N}}{|u(x_{j},T)-U(x_{j},T)|},\nonumber 
\end{eqnarray}
respectively. We set $x_{j}=\frac{j}{\mathcal{N}}$ and $\mathcal{N}=100$
in the computations. The spatial and temporal rate of convergence
of the method are computed by
\begin{eqnarray*}
 &  & rate_{N_{i}}=\frac{\log{\frac{E(N_{i},\tau)}{E(N_{i-1},\tau)}}}{\log{\frac{N_{i-1}}{N_{i}}}},\quad rate_{\tau_{i}}=\frac{\log{\frac{E(N,\tau_{i})}{E(N,\tau_{i-1})}}}{\log{\frac{\tau_{i}}{\tau_{i-1}}}},
\end{eqnarray*}
respectively, where $E(N,\tau)$ indicates the error with a basis
of dimension $N$ and time-step length $\tau$\@. However, we will
use the logarithmic scale plots to show the the method has a spectral
accuracy in space. The computations were performed by using Maple
18 on a Lenovo laptop running Windows 8.1 platform with a Core i3
1.90 GHz CPU and 4 Gb memory. 
\begin{example}
\label{exa:1}We consider the problem (\ref{eq:main}) with the homogeneous
boundary and initial conditions with $\kappa_{1}=0.1$ and $\kappa_{2}=2$.
The source term is such that the exact solution is $u=x^{2}(1-x)\sin{t}$.
Table \ref{tab:example1} shows the $L^{2}$ and $L^{\infty}$ errors
at time $t=1$ for different $N$ and $\tau=1/M$ with CPU time for
some fractional orders. From the table, we can see the convergence
of the proposed method. 
\begin{flushleft}
\begin{table}[h]
\caption{\label{tab:example1}$L_{2}$ and $L_{\infty}$errors for Example
\ref{exa:1} at $T=1$ with CPU time (s).}
\noindent\resizebox{\textwidth}{!}{%
\raggedright{}%
\begin{tabular}{ccccccccccc}
\hline 
 &  & \multicolumn{3}{c}{$\alpha=0.25$} & \multicolumn{3}{c}{$\alpha=0.5$} & \multicolumn{3}{c}{$\alpha=0.75$}\tabularnewline
\hline 
$M$ & $N$ & $L_{\infty}^{T}$ & $L_{2}^{T}$ & time & $L_{\infty}^{T}$ & $L_{2}^{T}$ & time & $L_{\infty}^{T}$ & $L_{2}^{T}$ & time\tabularnewline
\hline 
10 & 4    & 3.46E-5 & 7.11E-5 & 0.391 & 1.22E-4 & 2.45E-4 & 0.44 & 3.20E-4 & 6.21E-4 & 0.31\tabularnewline
20 & 6     & 1.45E-5 & 3.08E-5 & 1.547 & 6.09E-5 & 1.28E-4 & 1.11 & 1.97E-4 & 4.11E-4 & 1.31\tabularnewline
40 & 8     & 4.66E-6 & 9.58E-6 & 4.079 & 2.27E-5 & 4.64E-5 & 3.59 & 8.72E-5 & 1.76E-4 & 3.72\tabularnewline
80 & 10     & 1.46E-6 & 2.92E-6 & 10.844 & 8.29E-6 & 1.65E-5 & 11.28 & 3.77E-5 & 7.45E-5 & 11.00\tabularnewline
120 & 12     & 7.34E-7 & 1.44E-6 & 24.11 & 4.56E-6 & 8.90E-6 & 22.38 & 2.31E-5 & 4.47E-5 & 23.25\tabularnewline
160 & 14    & 4.49E-7 & 8.65E-7 & 38.922 & 2.98E-6 & 5.72E-6 & 41.56 & 1.62E-5 & 3.10E-5 & 40.75\tabularnewline
\hline 
\end{tabular}}
\end{table}
\par\end{flushleft}

\end{example}

\begin{figure}
\begin{centering}
\includegraphics[scale=0.4]{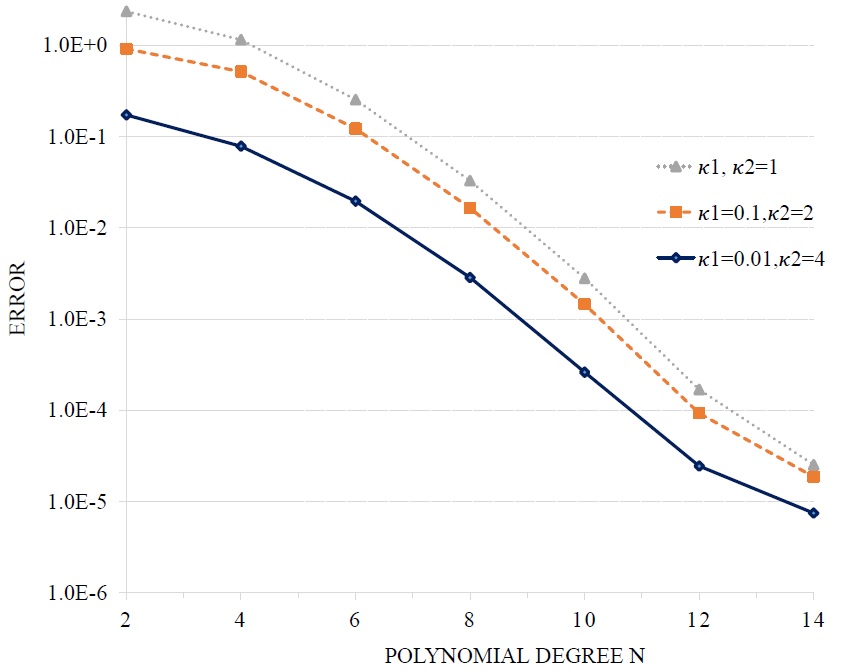}
\par\end{centering}
\caption{\label{fig:The-spectral-convergence}The $H^{1}$-error for some advection
and dispersion coefficients for Example \ref{exa:2}.}
\end{figure}

\begin{example}
\label{exa:2}For the problem (\ref{eq:main}), let $\alpha=0.5$,
the initial condition $g(x)=\sin{(\pi x)}$ with homogeneous boundary
conditions and the exact solution $u=\sin(\pi x)\exp{\left(-t^{2}\right)}.$
Fig. \ref{fig:The-spectral-convergence} illustrates the convergence
in space in $H^{1}-$norm for three cases of advection and dispersion
coefficients with $\tau=1/400$ at $t=1$\@. It is seen the logarithmic
scaled error behaves almost linearly versus the polynomial degree,
i.e., the so-called spectral accuracy. 
\end{example}

\begin{example}
\label{exa:3}Consider the problem (\ref{eq:main}) with $\kappa_{1}=0.2$
and $\kappa_{2}=1.5$ with the exact solution $u=x^{4}(1-x)^{2}t^{2}$.
Table \ref{tab:Ex3} shows the error results at $t=1$ for the fractional
orders $\alpha=0.25$, $\alpha=0.50$ and $\alpha=0.75$\@. 
\begin{flushleft}
\begin{table}
\caption{\label{tab:Ex3}Errors for Example \ref{exa:3} at $t=1$ with CPU
time (s). }
\noindent\resizebox{\textwidth}{!}{
\begin{tabular}{ccccccccccc}
\hline 
 &  & \multicolumn{3}{c}{$\alpha=0.25$} & \multicolumn{3}{c}{$\alpha=0.5$} & \multicolumn{3}{c}{$\alpha=0.75$}\tabularnewline
\hline 
$M$ & $N$ & $L_{\infty}^{T}$ & $L_{2}^{T}$ & time & $L_{\infty}^{T}$ & $L_{2}^{T}$ & time & $L_{\infty}^{T}$ & $L_{2}^{T}$ & time\tabularnewline
40 & 4 & 9.49E-2 & 2.14E-1 & 2.45 & 8.57E-2 & 1.90E-1 & 1.98 & 7.62E-2 & 1.66E-1 & 2.47\tabularnewline
80 & 6 & 3.01E-6 & 7.02E-6 & 5.94 & 1.86E-5 & 4.36E-5 & 5.78 & 9.56E-5 & 2.25E-4 & 5.58\tabularnewline
160 & 8 & 1.35E-7 & 2.65E-7 & 17.66 & 9.26E-7 & 1.83E-6 & 17.00 & 5.38E-6 & 1.08E-5 & 17.37\tabularnewline
320 & 10 & 4.22E-8 & 8.13E-8 & 53.88 & 3.35E-7 & 6.46E-7 & 49.27 & 2.27E-6 & 4.39E-6 & 52.86\tabularnewline
\hline 
\end{tabular}}
\end{table}
\par\end{flushleft}

\end{example}

\begin{example}
\label{exa:Ex4}Consider the problem (\ref{eq:main}) with $\kappa_{1}=0.1$,
$\kappa_{2}=2$ and the exact solution $u=x\cos(\frac{\pi}{2}x)\exp(-t)$.
The rate of convergence in time is reported in Table \ref{tab:TemporalOrderU=00003DxCos(..}
for $N=14$ at $t=1$. Also Fig. \ref{figEx4} illustrates the convergence
in space for $\alpha=0.5$ and $\tau=1/100$ at $t=1$ in terms of
$H^{1}-$norm. 

The examples confirm the theoretical results of convergence of the
time discretization (\ref{eq:TempOrder2-alpha}) and spectral discretization
(\ref{eq:SpatialRate}).

\begin{table}
\caption{\label{tab:TemporalOrderU=00003DxCos(..}Experimental rate of convergence
in time for Example \ref{exa:Ex4} at $t=1$ with $N=14$. }
\noindent\resizebox{\textwidth}{!}{
\begin{tabular}{ccccccccccccc}
\hline 
 & \multicolumn{4}{c}{$\alpha=0.25$} & \multicolumn{4}{c}{$\alpha=0.5$} & \multicolumn{4}{c}{$\alpha=0.75$}\tabularnewline
\hline 
$M$ & $L_{\infty}^{T}$ & rate & $L_{2}^{T}$ & rate & $L_{\infty}^{T}$ & rate & $L_{2}^{T}$ & rate & $L_{\infty}^{T}$ & rate & $L_{2}^{T}$ & rate\tabularnewline
\hline 
25 & 1.68E-5 &  & 3.51E-5 &  & 7.87E-5 &  & 1.64E-4 &  & 3.04E-4 &  & 6.34E-4 & \tabularnewline
50 & 5.12E-6 & 1.715 & 1.07E-5 & 1.715 & 2.79E-5 & 1.498 & 5.82E-5 & 1.497 & 1.27E-4 & 1.253 & 2.66E-4 & 1.252\tabularnewline
100 & 1.55E-6 & 1.720 & 3.24E-6 & 1.720 & 9.87E-6 & 1.498 & 2.06E-5 & 1.498 & 5.35E-5 & 1.252 & 1.12E-4 & 1.251\tabularnewline
200 & 4.71E-7 & 1.724 & 9.82E-7 & 1.723 & 3.50E-6 & 1.498 & 7.29E-6 & 1.498 & 2.25E-5 & 1.251 & 4.70E-5 & 1.251\tabularnewline
400 & 1.42E-7 & 1.724 & 2.98E-7 & 1.722 & 1.24E-6 & 1.498 & 2.58E-6 & 1.498 & 9.45E-6 & 1.250 & 1.97E-5 & 1.254\tabularnewline
\hline 
\end{tabular}}
\end{table}
\begin{figure}
\centering{}\includegraphics[scale=0.4]{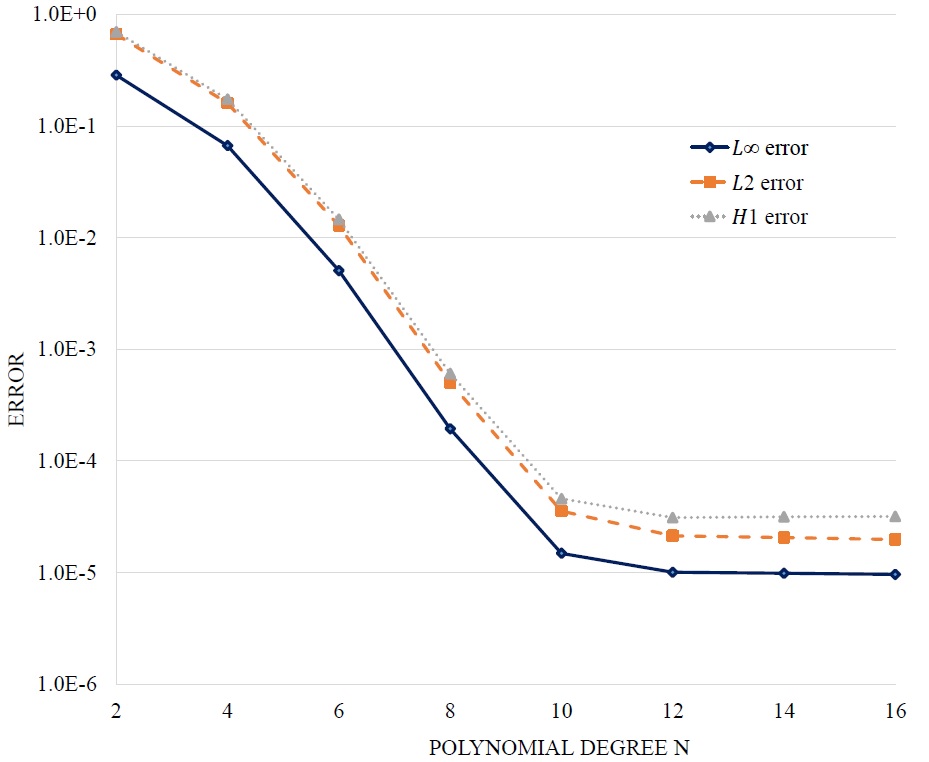}\caption{\label{figEx4}Convergence of the method for Example \ref{exa:Ex4}. }
\end{figure}
\end{example}

\section{\label{sec:Con}Conclusion}

In this paper, we stated a new property for derivatives of Bernstein
basis\@. Using this property, we obtained exact banded operational
matrices for derivatives of Bernstein basis\@. Since the basis transformation
may be ill-conditioned, the first advantage of our work in comparison
with the existing works is that we did not use any basis transformation
for the derivation of operational matrices. The second is that the
derived matrices are banded, so less computational effort is required
for a desired accuracy with less round-off errors. We also proposed
a numerical method utilizing the Petrov-Galerkin method for the time-fractional
advection dispersion equation on bounded domains based on the operational
matrices leading to banded linear systems. Moreover, we derived the
matrix formulation of the method and showed that the resulting linear
system has a unique solution. We also discussed the error analysis.
Then providing some numerical experiments, it is seen that the method
is efficient, accurate and simple to implement for solving time-fractional
advection dispersion equations on bounded domains.

\section*{Acknowledgment}

The authors would like to thank the anonymous reviewers for their
valuable suggestions and comments.

\end{document}